\documentstyle[12pt]{article}

\begin{document}
\topmargin-.75cm
\baselineskip16.125pt
\parskip.5\baselineskip

\renewcommand{\thesection}{\ifnum\value{section}=0 0\else
\arabic{section}\fi}

\newtheorem{definition}{Definition $\!\!$}[section]
\newtheorem{prop}[definition]{Proposition $\!\!$}
\newtheorem{lem}[definition]{Lemma $\!\!$}
\newtheorem{corollary}[definition]{Corollary $\!\!$}
\newtheorem{theorem}[definition]{Theorem $\!\!$}
\newtheorem{example}[definition]{Example $\!\!$}
\newtheorem{remark}[definition]{Remark $\!\!$}

\newcommand{\nc}[2]{\newcommand{#1}{#2}}
\newcommand{\rnc}[2]{\renewcommand{#1}{#2}}
\nc{\bpr}{\begin{prop}}
\nc{\bth}{\begin{theorem}}
\nc{\ble}{\begin{lem}}
\nc{\bco}{\begin{corollary}}
\nc{\bre}{\begin{remark}}
\nc{\bex}{\begin{example}}
\nc{\bde}{\begin{definition}}
\nc{\ede}{\end{definition}}
\nc{\epr}{\end{prop}}
\nc{\ethe}{\end{theorem}}
\nc{\ele}{\end{lem}}
\nc{\eco}{\end{corollary}}
\nc{\ere}{\end{remark}}
\nc{\eex}{\end{example}}
\nc{\epf}{\hfill\mbox{$\Box$}}
\nc{\ot}{\otimes}
\nc{\LAblp}{\mbox{\LARGE\boldmath$($}}
\nc{\LAbrp}{\mbox{\LARGE\boldmath$)$}}
\nc{\blp}{\mbox{\boldmath$($}}
\nc{\brp}{\mbox{\boldmath$)$}}
\nc{\LAlp}{\mbox{\LARGE $($}}
\nc{\LArp}{\mbox{\LARGE $)$}}
\nc{\Llp}{\mbox{\Large $($}}
\nc{\Lrp}{\mbox{\Large $)$}}
\nc{\llp}{\mbox{\large $($}}
\nc{\lrp}{\mbox{\large $)$}}
\nc{\lbc}{\mbox{\Large\boldmath$,$}}
\nc{\lc}{\mbox{\Large$,$}}
\nc{\Lall}{\mbox{\Large$\forall$}}
\nc{\bc}{\mbox{\boldmath$,$}}

\nc{\Section}{\setcounter{definition}{0}\section}
\renewcommand{\theequation}{\thesection.\arabic{equation}}

\newcounter{c}
\renewcommand{\[}{\setcounter{c}{1}$$}
\newcommand{\etyk}[1]{\vspace{-7.4mm}$$\begin{equation}\label{#1}
\addtocounter{c}{1}}
\renewcommand{\]}{\ifnum \value{c}=1 $$\else \end{equation}\fi}

\newcommand{\dowod}{\noindent{\bf Proof:} }
\newcommand{\Sp}{{\rm Sp}\,}
\newcommand{\Mor}{\mbox{$\rm Mor$}}
\newcommand{\skrA}{{\cal A}}
\newcommand{\Phase}{\mbox{$\rm Phase\,$}}
\newcommand{\id}{{\rm id}}
\newcommand{\diag}{{\rm diag}}
\newcommand{\inv}{{\rm inv}}
\newcommand{\ad}{{\rm ad}}
\newcommand{\poi}{{\rm pt}}
\newcommand{\Dim}{{\rm dim}\,}
\newcommand{\Ker}{{\rm ker}\,}
\newcommand{\Mat}{{\rm Mat}\,}
\newcommand{\Rep}{{\rm Rep}\,}
\newcommand{\Fun}{{\rm Fun}\,}
\newcommand{\Tr}{{\rm Tr}\,}
\newcommand{\supp}{\mbox{$\rm supp$}}
\newcommand{\half}{\frac{1}{2}}

\newcommand{\skrF}{{\cal F}}
\newcommand{\skrD}{{\cal D}}
\newcommand{\skrC}{{\cal C}}

\newcommand{\ttimes}{\mbox{$\hspace{.5mm}\bigcirc\hspace{-4.9mm}
\perp\hspace{1mm}$}}
\newcommand{\Ttimes}{\mbox{$\hspace{.5mm}\bigcirc\hspace{-3.7mm}
\raisebox{-.7mm}{$\top$}\hspace{1mm}$}}
\newcommand{\Cstar}{$^*$-}

\newcommand{\bbr}{{\bf R}}
\newcommand{\bbz}{{\bf Z}}
\newcommand{\Ci}{C_{\infty}}
\newcommand{\Cb}{C_{b}}
\newcommand{\fa}{\forall}
\newcommand{\te}{\exists}
\newcommand{\rrr}{right regular representation}
\newcommand{\wrt}{with respect to}

\newcommand{\qg}{quantum group}
\newcommand{\qgs}{quantum groups}
\newcommand{\cs}{classical space}
\newcommand{\qs}{quantum space}
\newcommand{\po}{Pontryagin}
\newcommand{\ch}{character}
\newcommand{\chs}{characters}
\newcommand{\hg}{\hat{G}}
\newcommand{\be}{\begin{equation}}
\newcommand{\ee}{\end{equation}}
\newcommand{\al}{\alpha}
\newcommand{\bet}{\beta}
\newcommand{\Lam}{$\Lambda$}
\newcommand{\La}{\Lambda}
\newcommand{\lam}{$\lambda$}
\newcommand{\la}{\lambda}
\newcommand{\ep}{\epsilon}
\newcommand{\Om}{\Omega}
\newcommand{\lr}{\longleftrightarrow}
\newcommand{\g}{\Gamma}
\newcommand{\gt}{{\widetilde{\Gamma}}}
\newcommand{\D}{\Delta}
\newcommand{\si}{\sigma}
\newcommand{\f}{\varphi}
\newcommand{\eps}{\varepsilon}
\newcommand{\w}{\omega}
\newcommand{\C}{{\bf C }}
\newcommand{\R}{{\bf R }}
\newcommand{\bea}{\begin{eqnarray}}
\newcommand{\eea}{\end{eqnarray}}
\newcommand{\ba}{\begin{array}}
\newcommand{\ea}{\end{array}}
\newcommand{\ber}{\begin{eqnarray}}
\newcommand{\ear}{\end{eqnarray}}
\newcommand{\U}{\Upsilon}
\newcommand{\der}{\mbox{d}}
\newcommand{\de}{\mbox{{\rm det}}}

\def\inbar{\,\vrule height1.5ex width.4pt depth0pt}
\def\IC{\relax\,\hbox{$\inbar\kern-.3em{\rm C}$}}
\def\otc{\otimes_{\IC}}
\def\ra{\rightarrow}
\def\ota{\otimes_ A}
\def\otza{\otimes_{ Z(A)}}
\def\otc{\otimes_{\IC}}
\def\h{\rho}
\def\x{\zeta}
\def\th{\theta}
\def\s{\sigma}
\def\t{\tau}
\def\st{\stackrel}
\title{\vspace*{-24mm}{\bf METRICS~AND~PAIRS~OF~LEFT~AND~RIGHT 
CONNECTIONS ON BIMODULES}}
\author{Ludwik D\c{a}browski\\
\normalsize SISSA, Via Beirut 2-4, Trieste, Italy.\\
\normalsize \sc e-mail: dabrow@sissa.it\vspace*{5mm} \\
Piotr M.~Hajac\thanks{
Partially supported by a visiting fellowship at 
the International Centre for Theoretical Physics in Trieste and by the KBN grant
\mbox{2 P301 020 07}. On~leave from:
Department of Mathematical Methods in Physics, 
Warsaw University, ul.~Ho\.{z}a 74, Warsaw, \mbox{00--682~Poland}.
{\sc e-mail: {\sc pmh@fuw.edu.pl}}
} \\
\normalsize ICTP, Strada Costiera 11, Trieste, Italy.\\
\normalsize \sc e-mail: pmh@ictp.trieste.it
\vspace*{5mm}\\ 
Giovanni Landi \\ 
\normalsize Dipartimento di Scienze Matematiche, Universit\'{a} di Trieste,
P.le Europa 1, Trieste, Italy\\
\normalsize 
and INFN, Sezione di Napoli, Napoli, Italy.
\\ 
\normalsize \sc e-mail: landi@univ.trieste.it\vspace*{5mm}\\
Pasquale Siniscalco\\
\normalsize SISSA, Via Beirut 2-4, Trieste, Italy.\\
\normalsize\sc e-mail: sinis@sissa.it \vspace*{5mm}
}
\date{\the\day\ February 1996}
\maketitle

\begin{abstract}
Properties of metrics and pairs consisting of left and right  connections are 
studied 
on the bimodules of differential 1-forms. Those bimodules are obtained from
the derivation based calculus of an algebra of matrix valued functions,
and an $SL\sb q(2,\IC)$-covariant calculus of the quantum plane plane at a 
generic $q$ and the cubic root of unity. It is shown that, in the aforementioned 
 examples, giving up the 
middle-linearity of metrics significantly enlarges the space of metrics.
A~metric compatibility condition for the pairs of left and right
connections is defined. Also, a compatibility condition between a left and right
connection is discussed. Consequences entailed by reducing to the centre of a
bimodule the domain of those conditions are investigated in detail.
Alternative ways of relating left and right connections are considered. 
\end{abstract}

\vfill
SISSA 26\slash 96\slash FM; QDSM-Trieste\slash 362; q-alg/9602035

\newpage\parindent5mm\baselineskip16.5pt\parskip.65\baselineskip
\Section*{Introduction}

Motivated to a great extent by the need to reconcile the geometric theory
of gravity
with the (noncommutative) operator algebraic theory of quantum physics,
there is considerable interest in generalising the formalism of General Relativity
to the realm of
Noncommutative Differential Geometry \cite{conbook}. 
In this paper, we study three concepts that are apparently
needed for such a generalisation: metric, linear connection and metric compatibility
condition.

We define a metric $g:E\times E\rightarrow A$ as a $\tau$-symmetric $A$-bilinear
pairing on an $A$-bimodule $E$, where $\tau$ is some generalised permutation.
Then we argue that
giving up the (often postulated) requirement that a metric
factor to a map defined on $E\otimes\sb A E$ one can obtain an essentially 
bigger space of metrics. In particular, we provide an example with 
an ample supply of $\tau$-symmetric metrics but where the requirement that a
$\tau$-symmetric metric $g$ descend to $E\ota E$ amounts to demanding that $g=0$
(see Proposition~\ref{nometric}).

Inspired by \cite{plane} on the one hand and by \cite{cq} on the other, 
we consider a pair of mutually compatible
connections on a bimodule. 
First connection of such a compatible pair
is a left connection in the sense that it satisfies the Leibniz rule with respect
to the left module structure. Similarly, the second connection is a right 
connection in the sense that it fulfills the Leibniz rule on the right. The 
compatibility condition is simply a requirement that the left and the right 
connection agree on the centre of a bimodule up to a 
bimodule isomorphism $\s$ 
(e.g., braiding). This bimodule isomorphism is, again, a generalised permutation.
Restricting the domain of
the aforementioned left-right compatibility condition to the centre of a bimodule
permits, at least in the considered examples, a significantly bigger space of 
solutions to this condition. (This seems desired at least
from the point of view of developing variational calculus on the space of 
connections.)

As to the compatibility between metrics and pairs of left and right
connections, we define a 2-parameter family of compatibility conditions, but then
restrict ourselves to the one that seems the most natural.

In the first section, we provide the general formalism and fix the notation.
Then we proceed to the first example, where $A$ is the algebra of matrix valued 
functions on a 
parallelizable manifold, and $E=A^1$ is the bimodule of 1-forms of the
derivation based differential calculus equipped with the 
pullback-of-permutation automorphism. 
Next we pass on to the quantum plane and the differential calculus 
with the braiding employed in~\cite{plane}.
First we consider the case of a generic $q$, and then the case of the cubic root
of unity. 
To obtain a non-zero $\tau$-symmetric 
metric on the quantum plane we have to `rescale' the
braiding $\s$ used in \cite{plane} by~$q^2$.
The thus obtained automorphism $\tau$ appears to be in a better agreement 
with the theory presented in~\cite{slw}.
(Even though in this case both automorphisms can be
considered over the same domain, they are different deformations of
the usual tensor product permutation.)

In what follows, Einstein's convention of summing over repeating indices is
assumed, and the unadorned tensor product stands for the tensor product over a field.

\parindent0pt\baselineskip16.125pt\parskip.5\baselineskip
\section{General Definitions}

Let $A$ be a unital associative algebra over a field $k$, and 
$E$ be a left and right projective
\mbox{$A$-b}imodule. We begin with a definition of a linear pairing 
$g:E\times E\to A$, which, for the sake of
simplicity (neglecting the nondegeneracy and reality conditions), 
we call a metric on $E$ (cf.~Section~8 and Section~9 in~\cite{dm2}).

\bde\label{g}
Let $\tau :E\ot E\to E\ot E$ be a bimodule automorphism. A linear map $g$
from  $E\ot E$ to $A$ is called a {\em $\tau$-symmetric metric} on $E$ 
(or simply {\em metric}, if no confusion arises) iff it is:\\
\hspace*{5mm}1)  bilinear over $A$, i.e.
$g(a\zeta , \rho b)=ag(\zeta , \rho )b\, ,\; 
\forall\, \x , \h \in E,\; a,b \in A\, ;$
\\
\hspace*{5mm}2) $\tau$-symmetric, i.e.~$g\circ\tau =g$.
\ede

Note that if $E$ is a central bimodule \cite{dm2},
that is, if the left and right multiplications by the elements of $Z(A)$ 
coincide on any element of $E$ 
(which is always the case in the examples considered in this article),
then any metric $g$ can be regarded as a map from $E \otza E$ to $A$.
We would like to emphasize here that, contrary to many other papers (e.g., see
(10) in~\cite{as}), we do not require $g$ to be well-defined on $E \ota E$.
As we show in our three examples, a requirement like this (which goes
under the name of middle-linearity) can be considered too restrictive 
(see Proposition~\ref{gmz}, Proposition~\ref{nometric} and Proposition~\ref{3ml}).
Giving up the middle-linearity condition allows us to get rid of those
restrictions.

Another structure on a bimodule $E$ that we wish to discuss 
is a pair of compatible left and right connections.

\bde\label{bi}
Let $(A\sp1,d)$ be a first order differential calculus on $A$, and
\[
\s : E \ota A^{1} \longrightarrow A^{1} \ota E
\]
be a bimodule isomorphism.
Also, let $\nabla^L$ be a left connection, i.e.~a linear map from $E$ to
$A^{1} \ota E$ satisfying the left Leibniz rule
\[
\nabla^L (a\zeta ) = da \ota\zeta + a \nabla^L\zeta\, ,\;
\forall\, a \in A,\,\zeta \in E, 
\]
and let $\nabla^R$ be a right connection, 
i.e.~a linear map from $E$ to $E \ota A^{1}$ 
fulfilling the right Leibniz rule
\[
\nabla^R (\zeta a) =  (\nabla^R\zeta)a + \zeta\ota da\, ,\;
\forall\, a\in A,\, \zeta\in E.
\]
A pair $(\nabla\sp L,\nabla\sp R)$ is called {\em $\s$-compatible} iff
\be\label{comp} 
\mbox{\large$\forall$}\,\zeta\in Z(E):\; 
\nabla^L \zeta = (\sigma\circ\nabla^R) \zeta\, ,
\ee
where $Z(E):=\{\zeta\in E \,|\; a\zeta = \zeta a,\,\fa\, a \in A \}$ 
is the centre of~$E$. 
\ede

Let us recall that in Section~8 of \cite{cq} a connection on a bimodule is
also defined as a pair consisting of a left and right connection. There,
however, instead of $\sigma$-compatibility condition~(\ref{comp}), the condition
of $\nabla\sp L$ being a right $A$-homomorphism and $\nabla\sp R$ being a 
left $A$-homomorphism is imposed. The latter condition, albeit it permits for
an interesting algebraic theory, cannot be directly transferred to the commutative
case $Z(E)=E$. On the other hand, defining a connection on a bimodule by requiring,
much as in Definition~3.2 in~\cite{glq}, that the $\s$-compatibility condition is 
fulfilled on the
whole bimodule $E$ rather than just its centre~$Z(E)$, automatically yields, for the 
appropriate $\s$ (i.e.~the usual tensor product flip), the
standard definition of a connection in the classical case, but entails essential
restrictions on the space of connections in noncommutative examples 
(see~(\ref{cmz}), (\ref{french}), (\ref{3f}), \cite{plane2} and 
Theorem~5.6 in~\cite{glq}; cf.~Lemma~1 in~\cite{dm2}).
If we demand that the equality \mbox{$\nabla\sp L=\s\circ\nabla\sp R$}
be satisfied on the whole bimodule~$E$, we can equivalently think of a pair 
$(\nabla\sp L,\nabla\sp R)$ as a left connection $\nabla\sp L$ fulfilling an
additional (right) Leibniz rule of the form
\[
\nabla^L (\zeta a) =  (\nabla^L\zeta)a + \s (\zeta\ota da)\, ,\;
\forall\, a\in A,\, \zeta\in E,
\]
(cf.~the Introduction in~\cite{ext}). In the classical differential geometry,
with the help of the tensor product flip, any left connection uniquely determines
the corresponding right connection, and vice-versa, without imposing any 
limitations on either of the connections. As we demonstrate in the next section,
this is precisely what happens with pairs of \mbox{$\s$-compatible} 
connections in the noncommutative
example of a `matrix geometry' (see Proposition~\ref{matrixbi}). A very similar
result (Proposition~\ref{3bi}) 
is obtained for a bimodule of differential 1-forms on the quantum plane
(see Proposition~2 in~\cite{sch}) at the cubic root of unity. (In fact, in all 
examples presented in this paper, we put $E=A\sp1$, so that the pairs of 
\mbox{$\s$-compatible} connections
studied here can be thought of as linear connections.)
When dealing with connections and metrics, it seems that in both cases we have
the same mechanism at work: solving constraints over just the commutative part allows
solutions to be parametrized by a whole noncommutative algebra, whereas solving
them over an entire noncommutative space renders the solutions parametrized by
the centre of an algebra. The examples of subsequent sections allow us
to explore this mechanism in some detail.

Let us now consider possible compatibility conditions between 
a pair of connections $(\nabla^{L}, \nabla^{R})$ and a metric $g$.
In order to do so, first
we must define two extensions of 
$g$ :
\[
\check{g} : A^1 \otimes\sb A E \otimes E \longrightarrow A^1  \;\mbox{\ and\ }\;
\hat{g} :  E \otimes E \otimes\sb A A^1\longrightarrow A^1 \, .
\]
It appears natural to choose:
\begin{equation}\check{g}(\alpha\otimes\sb A\zeta ,\rho) 
=\alpha g(\zeta ,\rho)\, ,\hspace{1cm}
\forall\alpha\in A^1,\;\zeta ,\rho\in E; 
\end{equation}\begin{equation}
\hat{g}(\zeta ,\rho\otimes\sb A \alpha)=g(\zeta ,\rho)\alpha\, ,\hspace{1cm}
\forall\alpha\in A^1,\;\zeta ,\rho\in E. 
\end{equation}

In principle, one can formulate the class of metric compatibility conditions
 by requiring the diagram

\begin{equation}\label{diagram}
\def\normalbaselines{\baselineskip30pt
\lineskip3pt \lineskiplimit3pt }
\def\mapright#1{\smash{
\mathop{
{-\!\!\!}-\!\!\!{-\!\!\!}{-\!\!\!}-\!\!\!{-\!\!\!}{-\!\!\!}
-\!\!\!{-\!\!\!}{-\!\!\!}\!\!\!{-\!\!\!}{-\!\!\!}-\!\!\!{-\!\!\!}-\!\!\!
\longrightarrow\!\!\!}
\limits^{#1}}}
\def\mapdown#1{\Bigg\downarrow
\rlap{$\vcenter{\hbox{$\scriptstyle#1$}}$}}
\matrix{E\otimes E&\mapright{\ f\sb L(t)\otimes id\,\oplus\, id\otimes f\sb R(s)}
&(\! A\sp1\!\otimes\sb{\! A}\! E\!\otimes\! E)\oplus 
(\! E\!\otimes\! E\!\otimes\sb{\! A}\! A\sp1)
\cr
\mapdown{\mbox{$g$}}&&\hidewidth\mapdown{\mbox{$\check{g}\oplus\hat{g}$}}\hidewidth
\cr
A&\stackrel{\mbox{$d\ \ $}}
{\!\!\!{-\!\!\!}-\!\!\!{-\!\!\!}-\!\!\!{-\!\!\!}\longrightarrow}\;\;\; A\sp1\;\;\;
\stackrel{\mbox{$\ \ +$}}
{\longleftarrow\!\!\!{-\!\!\!}-\!\!\!{-\!\!\!}-\!\!\!{-\!\!\!}}
\hidewidth
&A\sp1\oplus A\sp1 , 
\cr
}
\end{equation}

where $f\sb L(t):=(1-t)\nabla\sp L+t(\s\circ\nabla\sp R)$, 
$f\sb R(s):=s(\s\sp{-1}\circ\nabla\sp L)+(1-s)\nabla\sp R$, 
$t,s\in k$, to commute. Here, however, in order to ensure that $f\sb L(t)$ and
$f\sb R(s)$ are connections, we settle for the particular case
$t=0=s$:

\bde\label{g-comp}
We say that a pair of connections {\em $(\nabla\sp L,\nabla\sp R)$
is compatible with $g$}~iff
\be\label{g-comp-f}
 dg(\zeta,\rho) = \check{g}(\nabla^{L}\zeta,\rho) +
\hat{g}(\zeta,\nabla^{R} \rho),\ \;\forall\,\zeta,\rho
\in E. 
\ee
\ede

As we show in Proposition~\ref{exte}, formula (\ref{g-comp-f}) is not always
sensitive to whether we consider it over $E$ or only over $Z(E)$. Observe that, 
if (\ref{comp}) is satisfied, then on the centre of a bimodule
 we have $f\sb L(t)=\nabla^{L}$ and
$f\sb R(s)=\nabla^{R}$ for any values of $t,s\in k$, and, consequently, all
metric compatibility conditions for a pair of $\s$-compatible connections
are equivalent when considered over~$Z(E)$. 
Finally, let us remark that
the metric compatibility condition for a pair of connections related by
$\nabla^{L}= \sigma \circ \nabla^{R}$ (on the whole bimodule) that is given by (1.21)
and (1.26) in \cite{plane} seems inappropriate for the cases when the metric is not
middle-linear.

Next, we apply the above definitions in some quantum geometric models.

\Section{Algebras of Matrix Valued Functions}

Let us choose
$
A=C^{\infty}(M)\otimes M_n(\IC),\; E=A\sp1=\mbox{Hom}\sb{Z(A)}(\mbox{Der} A, A),\;
$\[
\tau(\alpha\otimes\beta )(X,Y)=(\alpha\otimes\beta )(Y,X),
\;\forall\, X,Y\in\mbox{Der} A,
\]
and $\s$ equal to $\tau$ factored to an automorphism of $A\sp1\ota A\sp1$. 
Here $M$ is a parallelizable manifold of dimension $m$, and the ground field of
$A$ is the field of complex numbers. Now, let 
$\{\theta\sp i\}\sb{i\in\{ 1,...,m+n\sp2-1\}}$ be the basis of $A\sp1$ as
defined in Section~3 of~\cite{ext}.
An important property of this basis is that 
\be\label{aa}
a\theta\sp i=\theta\sp ia,\;\;\;\forall\, a\in A,\, i\in\{1,...,m+n\sp2-1\}\, ,
\ee 
\be\label{aa2}
\tau(\theta\sp i\otimes\theta\sp j)=\theta\sp j\otimes\theta\sp i,
\;\;\;\forall\, a\in A,\, i,j\in\{1,...,m+n\sp2-1\}\, .
\ee 
In this setting, one can immediately verify
the following claim (cf.~Section~9 in~\cite{dm2} and p.5861 in~\cite{ext}):

\bpr\label{gmz}
Let $g\sp{ij}$ denote $g(\theta\sp i\otimes\theta\sp j)$, where 
$i,j\in\{1,...,m+n\sp2-1\}$. The map 
$
\psi :g\mapsto\left( g\sp{ij}\right)
$ 
provides a one-to-one correspondence between the metrics (the 
middle-linear metrics) on $A\sp1$ and
the symmetric matrices of $M\sb{m+n\sp2-1}(A)$ (the symmetric matrices of 
$M\sb{m+n\sp2-1}(Z(A))$ respectively). 
\epr

In the same basis, let us define the Christoffel symbols of $\nabla\sp L$
and $\nabla\sp R$ by
\begin{eqnarray}\label{tof}
\nabla\sp L\theta\sp i &=& \theta\sp j\ota\theta\sp k\g\sp i\sb{jk}\, ,
\nonumber\\
\nabla\sp R\theta\sp i &=& \theta\sp j\ota\theta\sp k\widetilde{\g}\sp i\sb{jk}\, .
\end{eqnarray}
Taking into account (\ref{aa}) and (\ref{aa2}) and noticing that
$\{\theta\sp i\}\sb{i\in\{ 1,...,m+n\sp2-1\}}$ is also a basis of the
\mbox{$Z(A)$-module} $Z(A\sp1)$, it is straightforward to prove:
 
\bpr\label{matrixbi}
 
A pair of connections $(\nabla\sp L,\nabla\sp R)$ is $\s$-compatible if and only
if its Christoffel symbols satisfy the equation
\be\label{mz}
\widetilde{\Gamma}^i_{kj} = \Gamma^i_{jk}\, .
\ee
Similarly, a pair of connections is $\s$-compatible on the whole bimodule
if and only if ((3.9) in~\cite{ext})
\be\label{cmz}
\widetilde{\Gamma}^i_{kj} = \Gamma^i_{jk}\in C\sp\infty(M)\, .
\ee
\epr

Thus the $\s$-compatibility condition (over $Z(E)$) allows $\nabla\sp L$ to
uniquely determine $\nabla\sp R$ and vice-versa. This is not unexpected since
(in this set-up) 
\be\label{azza}
A^1 = A Z(A^1) = Z(A^1)A
\ee
and $\nabla\sp L$ and $\nabla\sp R$ satisfy the left and right
Leibniz rule respectively. On the other hand, as we shall see in Section~4 
(Remark~\ref{azzar})
$\nabla\sp L$ and $\nabla\sp R$ can mutually determine each other even if 
(\ref{azza}) is not satisfied. 

Concerning the metric compatibility of $(\nabla^L,\nabla^R)$, we can again take
an advantage of (\ref{aa}) to show that (\ref{g-comp-f}) is equivalent to
\be\label{dg}
 dg^{ij}=(\g^i_{kl} g^{lj} + g^{il}\widetilde{\g}^j_{lk})\theta^k\, ,
\;\;\;\forall\, i,j\in\{1,...,m+n\sp2-1\}\, . 
\ee
To end this section, let us observe that, if $(\nabla\sp L,\nabla\sp R)$ is
$\s$-compatible, then (\ref{dg}) coincides with (3.13) in~\cite{ext}. One should
bear in mind, however, that the latter has been obtained from a different
starting point ((1.9)~in~\cite{ext}) and only for middle-linear metrics. 

\Section{Generic Quantum Plane}

The next example that we study regards a bimodule of differential 1-forms ($E=A\sp1$)
on the quantum plane. We choose as our space of differential 1-forms the grade one
of the differential algebra $\Omega(A)=A\oplus A\sp1\oplus A\sp2$ 
(e.g., see~\cite{plane}) 
that is given by the generators $1,x,y,\xi,\eta$, where 
$\xi=\th\sp1=dx,\,\eta=\th\sp2=dy$, and relations
\bea
&& x y = q y x ~, \nonumber \\
&& x \xi = q^2 \xi x ~, ~~~ x \eta = q \eta x + (q^2-1) \xi y~, ~~~
 y \xi = q \xi y ~, ~~~  y \eta = q^2 \eta y ~, \nonumber \\
&& \eta \xi + q \xi \eta = 0~, ~~~ \xi^2 = 0~, ~~~\eta^2 = 0~.
\label{comm}
\eea
For our bimodule automorphism $\s$ it is  natural (see the paragraph between
(2.11) and (2.12) in~\cite{plane}) to take the map 
defined by 
\be\label{planeflip}
\begin{array}{ll}
\s(\xi \ota \xi) = q^{-2} \xi \ota \xi~, &
\s(\xi \ota \eta) = q^{-1} \eta \ota \xi~, \\ 
\s(\eta \ota \xi) = q^{-1} \xi \ota \eta - (1-q^{-2}) \eta \ota \xi ~, &
\s(\eta \ota \eta) = q^{-2} \eta \ota \eta~. \\ 
\end{array}  
\ee
First we consider  the case of a generic $q$.
Then the centre of $A$ is $\IC$ and, as can be checked by a direct computation,  
the centre of $A^1$ is zero. If we choose $\tau$ equal to $\s$ (modulo 
the tensor product over $A$, as was done in the previous
section), we can immediately see that, unless $q=1$, there exists no non-zero
metric~\cite{plane}. To remedy this problem, we `rescale'~$\s$ by~$q\sp2$. 
More precisely, we put  
\be\label{tau}
\begin{array}{ll}
\tau(\xi \otimes \xi) = \xi \otimes \xi~, &
\tau(\xi \otimes \eta) = q\eta \otimes \xi~, \\ 
\tau(\eta \otimes \xi) = q\xi \otimes \eta - (q^{2}-1) \eta \otimes \xi ~, &
\tau(\eta \otimes \eta) = \eta \otimes \eta~. \\ 
\end{array}  
\ee
It turns out that this $\tau$ is quite natural from
the point of view of \cite{slw} --- it factors to an automorphism of
$A\sp1\ota A\sp1$ and preserves $\th\ota\th$, where
$\th =x\eta-qy\xi$ is the only (up to a multiplication by a complex number)
left and right \mbox{$SL\sb q(2,\IC)$-co}invariant 1-form 
(see Section~2 in~\cite{plane}). 
(Recall that if $A$ is a Hopf
algebra, then there exists a unique bimodule homomorphism $\tau$ such that
$\tau(\alpha\sb L\ota\alpha\sb R)=\alpha\sb R\ota\alpha\sb L$ for any left
coinvariant 1-form $\alpha\sb L$ and right coinvariant 1-form $\alpha\sb R$
--- see Proposition~3.1 in~\cite{slw}.)

It is obvious that with $\tau$ specified as above, the only constraints that
the coefficients of a metric $g$ have to satisfy is
\be\label{symm}
 g(\xi, \eta) = qg( \eta, \xi ) ~. 
\ee
On the other hand, it can be computed that the middle-linearity of $g$ is equivalent 
to the following equations:
\bea\label{mlqp}
xg(\xi, \xi ) & = & q^4g(\xi, \xi )x~, \nonumber \\
yg(\xi, \xi ) & = & q^2g(\xi, \xi )y~, \nonumber \\
xg(\xi, \eta) & = & q^3 g(\xi, \eta)x + 
 q^2 (q^{2}-1)g(\xi, \xi )y~, \nonumber \\ 
yg(\xi,  \eta) & = & q^3 g(\xi,  \eta)y~, \nonumber \\
xg(\eta, \xi ) & = & q^3 g( \eta, \xi )x + q(q^{2}-1)g(\xi , 
\xi )y~, \nonumber \\
yg( \eta, \xi ) & = & q^3 g( \eta, \xi )y~,  \nonumber \\
xg(  \eta, \eta) & = & q^{2}g( \eta, \eta)x +q^{2}(q^2-1)g( \xi, \eta)y+ 
q(q^2-1)g( \eta, \xi)y~, \nonumber\\
yg( \eta, \eta) & =& q^4g( \eta, \eta)y~.
\eea
Notice that, due to the commutation relation $xy=qyx$,
 as long as $q$ is generic and no negative powers of $x$ and $y$ 
are allowed, there is no non-zero solution to (\ref{mlqp}) 
(look at the equations with $y$). Consequently, we obtain:

\bpr\label{nometric}
If $q$ is not a root of unity, there is no ($\tau$-symmetric) middle-linear metric 
on $A\sp1$.  
\epr

\bre\label{neg}\em
If we admit negative powers of $x$ and $y$, the solutions of (\ref{mlqp})
form the following three-parameter family:
\bea
g(\xi, \xi ) & = & ax^{-2}y^4~, \nonumber \\
g(\xi, \eta) & = & qx^{-3}\left(by^{3} + q^{3}ay^5\right)~,
\nonumber \\ 
g( \eta, \xi ) & = & x^{-3}\left(by^{3} + q^{3}ay^5\right)~,
 \nonumber \\
g( \eta, \eta) & = & x^{-4}\left(cy^{2} +q^3(q^2 +1)by^4
+ q^8ay^6\right)\, ,
\eea
where $a$, $b$, $c$ are complex parameters.
\hfill{$\Diamond$}\ere

Since the centre of $A\sp1$ is zero, we can immediately conclude that

\bpr\label{genericbi}
If $q$ is not a root of unity, any pair $(\nabla^L,\nabla^R)$ of 
\mbox{$\s$-compatible} connections on
$A^1$ is a pair of independent and unrestricted left and right connections. 
\epr

On the other hand, there exists only a one-parameter family of solutions of the 
$\s$-compatibility condition considered over the whole $A\sp1$ 
(see~(2.13) in~\cite{plane}). 
Defining Christoffel symbols \linebreak
$\{ F\sp i\sb{jk},\,\widetilde{F}\sp i\sb{jk}\}\sb{i,j,k\in\{ 1,2\}}$ of such a
compatible pair of connections as in (\ref{tof}), we can write the aforementioned
solutions in the following way:  \ 
\bea\label{french}
&& F^1_{11} = \nu q x y^2~, ~~~F^1_{12} = -\nu q^3 x^2 y~, ~~~
F^1_{21} = -\nu q^2 x^2 y~, ~~~ F^1_{22} = \nu q^5 x^3~, 
\nonumber \\ 
&& F^2_{11} = \nu q^3 y^3~, ~~~F^2_{12} = -\nu q^4 x y^2~, ~~~
F^2_{21} = -\nu q^3 x y^2~, ~~~ F^2_{22} = \nu q^5 x^2 y ~,
\eea
where $\nu$ is a complex parameter. The Christoffel symbols 
$\{\widetilde{F}\sp i\sb{jk}\}\sb{i,j,k\in\{ 1,2\}}$ can be expressed in a
similar fashion.

As to the metric compatibility condition, formula (\ref{g-comp-f}) reads
\be\label{compa}
a_{i} \llp dg(\theta^i,\theta^j)-\theta^k g(\theta^l \g ^i_{kl}, \th^j)
-g(\th^i,\th^k)\th^l \gt^j_{kl}\lrp \widetilde{a}_j =0,\;\;\;
\forall\, a_{i},\widetilde{a}_j\in A,\, i,j\in\{ 1,2\}.
\ee
Clearly, (\ref{compa}) is satisfied if and only if the expression in the large
 parentheses vanishes for any $i$ and $j$.

\Section{Quantum Plane at the Cubic Root of Unity}

The setting of this section is identical with the setting of the previous one
except that now we take $q=e\sp{\frac{2\pi i}{3}}$ rather than generic $q$.
(For the sake of simplicity, we call $e\sp{\frac{2\pi i}{3}}$ {\it the} cubic
root of unity.)
Long but rather straightforward reasoning enables one to prove
the following lemma:

\ble\label{compute}
Let $q$ be the cubic root of unity. Then
\bea\label{za} \ &&
Z(A) = \{a_{ij} x^{3i} y^{3j}\, |\; a_{ij} \in \IC\},
\\ \label{dza} \ &&
dZ(A)=0,
\\ \label{aza} \ &&
\mbox{\large$\forall$}\, a\in Z(A),\, \alpha\in A\sp1:\; a\alpha=\alpha a,
\\ \label{za1} \ && 
Z(A^1)=\{ c\sb i\th\sp i\in A\sp1\, |\;
c_1=axy-bxy^3,\, c_2=bx^2y^2,\, a,b\in Z(A)\},
\\ \label{za2} \ &&
c\sb i\th\sp i=\th\sp j\tilde{c}\sb j,\;\mbox{where}\;\; c_1=axy-bxy^3,\,
c_2=bx^2y^2,\,
\nonumber \\ &&
\tilde{c}_1=axy-qbxy^3,\,\tilde{c}_2=c_2,\, a,b\in Z(A)
\;\mbox{(see~(\ref{cf2}))}.
\eea
\ele

Changing $q$ from generic to $q=e\sp{\frac{2\pi i}{3}}$ entails no consequence
as far as the (general $\tau$-symmetric) metric is concerned. However, regarding
middle-linear metrics, with the
help of commutation formulas provided in the appendix, one can prove:

\bpr\label{3ml}
If $q$ is the cubic root of unity and $g$ is middle-linear (but not
necessarily $\tau$-symmetric), 
then (\ref{mlqp}) is equivalent to:
\bea\label{g3z}
g(\xi, \xi ) & = & x^{3}Zxy~, \nonumber \\
g(\xi, \eta) & = & qx^{3}Zy^{2} + x^{3}Y~, \nonumber \\
g( \eta, \xi ) & = & x^{3}Zy^{2} + x^{3}W~, \nonumber \\
g( \eta, \eta) & = & Ux^{2}y^{2} +(qY + W)x^{2}y + 
q^{2}Zx^{2}y^{3}~,
\eea
where $Z$, $Y$, $W$, $U$ are arbitrary elements of $Z(A)$. 
Furthermore, $g$ is $\tau$-symmetric if and only if $Y=qW$.
\epr

Thus, much as in Proposition~\ref{gmz}, the space of middle linear metrics is 
three-dimensional over~$Z(A)$. This is not unexpected, if one remembers that
the quantum plane at the $n$-th root of unity is nothing but 
$\IC [x,y]\otimes M\sb n(\IC)$ (cf.~Section~IV.D.15 of~\cite{hw}).

Our next step is to determine the space of pairs of $\s$-compatible connections.
In order to make our reasoning more transparent, we introduce {\em formal}
inverses of $x$ and $y$. (It is simply more convenient, for instance, to write
$\widetilde\g =xy\g x\sp{-1}y\sp{-1}$ as the solution of the
equation $\widetilde\g xy=xy\g$ rather than consider $\widetilde\g$ and
$\g$ as power series in $x$ and $y$ and then express  
the complex coefficients of $\widetilde\g$
in terms of the complex coefficients of $\g$. However, we neither need nor assume
the existence of $x\sp{-1}$ and $y\sp{-1}$ in our algebra.) Treating 
$\{\th\sp i\ota\th\sp j\}\sb{i,j\in\{1,2\}}$ as a basis of the right $A$-module
$A\sp1\ota A\sp1$ and taking advantage of Lemma~\ref{compute}, one can carry out
lengthy but straightforward calculations that show that the $\s$-compatibility
condition (\ref{comp}) is equivalent to:

\vspace*{-3mm}\baselineskip12pt\begin{eqnarray}
&& \gt ^1_{11} =
q^2 x y \g ^1_{11} y^{-1}x^{-1}
+ (1\! -\! q) y^2 \g ^1_{12} y^{-1} x^{-1}
+ (q^2 \! -\! 1) y^2 \g ^1_{21} y^{-1} x^{-1}~, \nonumber \\
&& ~ \nonumber \\
&& \gt ^1_{12} =
(q^2\! -\! 1) x y \g ^1_{12} y^{-1} x^{-1}
+ q x y \g ^1_{21} y^{-1}x^{-1} 
+(q\! -\! q^2 ) y^2 \g ^1_{22}  y^{-1} x^{-1}~, \nonumber \\
&& ~ \nonumber \\
&& \gt ^1_{21} = q x y \g ^1_{12} y^{-1}x^{-1}
+ (1\! -\! q) y^2 \g ^1_{22} y^{-1} x^{-1}~, \nonumber \\
&&~ \nonumber \\
&& \gt ^1_{22} = q^2 x y \g ^1_{22} y^{-1}x^{-1}~, \nonumber \\
&&~ \nonumber \\
&& \gt ^2_{11} = x y \g ^1_{11} x^{-2} - 
x \g ^1_{11} yx^{-2}
+ (q\! -\! 1)  y\g ^1_{12} yx^{-2}
+ (q\! -\! q^2 )  y^2 \g ^1_{12} x^{-2} \nonumber \\
&& ~~~~~~~~ + (1\! -\! q^2 ) y\g ^1_{21} yx^{-2}
+ (1\! -\! q ) y^2 \g ^1_{21} x^{-2}
+ q^2 x^2 y^2 \g ^2_{11} y^{-2}x^{-2} \nonumber\\
&& ~~~~~~~~ +(q\! -\! 1) x \g ^2_{12} yx^{-2}
+ (1\! -\! q^2 )x \g ^2_{21} yx^{-2}
+ 3 y\g ^2_{22} yx^{-2}~, \nonumber \\
&&~\nonumber \\
 && \gt ^2_{12} =
(1\! -\! q^2 )x \g ^1_{12} yx^{-2}
+(q^2 \! -\! 1)y^2 \g ^1_{22} x^{-2}
+(1\! -\! q )  x y \g ^1_{12} x^{-2}\nonumber \\
&& ~~~~~~~~+ q^2  x y \g ^1_{21} x^{-2}
- q x \g ^1_{21} yx^{-2}
+ (q\! -\! 1 ) y\g ^1_{22} yx^{-2}\nonumber \\
&& ~~~~~~~~
+ (q^2\! -\! 1 ) x^2 y^2 \g ^2_{12} y^{-2}x^{-2}
+q x^2 y^2 \g ^2_{21} y^{-2}x^{-2}
+(q^2\! -\! q ) x \g ^2_{22} yx^{-2}~, \nonumber \\
&&~ \nonumber \\
&&\gt ^2_{21} =
  q^2 x y \g ^1_{12} x^{-2}
- q x  \g ^1_{12} yx^{-2}
+ (1\! -\! q^2 ) y\g ^1_{22} yx^{-2}\nonumber \\
&& ~~~~~~~~+ (q \! -\! q^2 )  y^2 \g ^1_{22} x^{-2}
+ q x^2 y^2 \g ^2_{12} y^{-2}x^{-2}
+ (q\! -\! 1) x \g ^2_{22} yx^{-2}~, \nonumber \\
&& ~ \nonumber \\
&& \gt ^2_{22} =  x y \g ^1_{22} x^{-2} - q x \g ^1_{22} y x^{-2}
+ q^2 x^2 y^2 \g ^2_{22} y^{-2}x^{-2}~,
\label{sol}
\end{eqnarray}\baselineskip16.125pt

where the Christoffel symbols are defined as in~(\ref{tof}). The above system of
equations allows one to determine uniquely $\nabla\sp R$ through $\nabla\sp L$,
but, as can be seen from the powers of $x$ and $y$, it cannot be done for an
arbitrary left connection $\nabla\sp L$. (The total power of $x$ and the total 
power of $y$ in each term of the right hand side of (\ref{sol}) have to be 
non-negative in order for (\ref{sol}) to make sense.) Since only total powers 
of $x$ turn negative in (\ref{sol}), it is convenient to think of an element of
$A$ as a polynomial in $x$ with coefficients in polynomials in $y$:
\be\label{gsplit}
\g ^i_{jk} = x\sp l\g ^i_{jkl}\, ,\;\;\; i,j,k\in\{ 1,2\}. 
\ee
To determine the necessary and
sufficient conditions  that $\{\g\sp i\sb{jk}\}\sb{i,j,k\in\{ 1,2\}}$ have to satisfy
in order to make (\ref{sol}) well-defined on the quantum plane, we substitute
(\ref{gsplit}) to (\ref{sol}) and conclude that the necessary and
sufficient conditions  for $\{\g\sp i\sb{jk}\}\sb{i,j,k\in\{ 1,2\}}$ are fully
given by:

\vspace*{-3mm}\baselineskip12pt\begin{eqnarray}
&&
(1\! -\! q)y^2\g^1_{120}+(q^2\! -\! 1)y^2\g^1_{210}=0~,\nonumber\\
&& ~ \nonumber \\
&&
(q\! -\! q^2 ) y^2\g^1_{220}=0~, \nonumber \\
&& ~ \nonumber \\
&&(1\! -\! q)y^2\g^1_{220}=0 ~, \nonumber \\
&&~ \nonumber \\
&&xy\g^1_{110}-x\g^1_{110}y
+ (q\! -\! 1)  y\g ^1_{120} y+ (q\! -\! 1)  yx\g ^1_{121} y
+ (q\! -\! q^2 )  y^2 \g ^1_{120}\nonumber \\
&&+ (q\! -\! q^2 )  y^2x \g ^1_{121}
+ (1\! -\! q^2 ) y\g ^1_{210} y+ (1\! -\! q^2 ) yx\g ^1_{211} y
+ (1\! -\! q ) y^2 \g ^1_{210}\nonumber\\
&&+ (1\! -\! q ) y^2 x\g ^1_{211} 
+(q\! -\! 1)x\g^2_{120} y+ (1\! -\! q^2 )x\g^2_{210} y
+ 3 y\g^2_{220} y+ 3 yx\g^2_{221} y=0~, \nonumber \\
&&~\nonumber \\
&&(1\! -\! q^2 )x \g ^1_{120} y
+(q^2 \! -\! 1) y^2 \g ^1_{220} +(q^2 \! -\! 1) y^2 x\g ^1_{221}
+(1\! -\! q )  x y \g ^1_{120}+ q^2  x y \g ^1_{210}~\nonumber \\
&& - q x \g ^1_{210} y
+ (q\! -\! 1 ) y\g ^1_{220} y+ (q\! -\! 1 ) yx\g ^1_{221} y
+(q^2\! -\! q) x \g ^2_{220} y = 0~, \nonumber \\
&&~ \nonumber \\
&&q^2 x y \g ^1_{120} - q x  \g ^1_{120} y
+ (1\! -\! q^2 ) y\g ^1_{220} y+ (1\! -\! q^2 ) yx\g ^1_{221} y\nonumber \\
&&+ (q \! -\! q^2 )  y^2 \g ^1_{220}+ (q \! -\! q^2 )  y^2 x\g ^1_{221} 
+ (q\! -\! 1) x \g ^2_{220} y=0~, \nonumber \\
&& ~ \nonumber \\
&&x y \g ^1_{220}  - q x \g ^1_{220} y=0~.
\label{sol2}
\end{eqnarray}\baselineskip16.125pt

Those equations are equivalent to: 
\bea\label{preprop}
&&
\g\sp1\sb{210}=\g\sp1\sb{120}=\g\sp1\sb{220}=\g\sp1\sb{221}=\g\sp2\sb{220}=0 
\nonumber \\ &&
(q-1)\g\sp2\sb{120} + (1-q^2)\g\sp2\sb{210}+ 3q^2y\g\sp2\sb{221}=0
\eea

Thus we have obtained:

{\samepage
\bpr\label{3bi}
The $\s$-compatibility condition (\ref{comp}) has a solution if and only if 
\bea\label{3bif}
&& \g ^1_{11} \in A~,~~~~~~~~~~
 \g ^1_{12} \in xA~,~~~~~~~~~~
 \g ^1_{21} \in xA~,~~~~~~~~~~ 
\g ^1_{22} \in x^2A~, \nonumber \\
&& \g ^2_{11} \in A~,~~~~~
 (q-1)\g ^2_{12}+(1-q\sp2)\g^2_{21}+
3q\sp2yx\sp{-1}\g^2_{22} \in xA~,~~~~~
\g ^2_{22} \in xA~.~~~~~~~~~~
\eea\nopagebreak
Moreover, if (\ref{3bif}) is satisfied, then the general solution of
(\ref{comp}) is given by (\ref{sol}).
\epr}

\bre\label{bire}\em
Expressing $\{\g ^i_{jk}\}\sb{i,j,k\in\{1,2\}}$ in terms of 
$\{\gt ^i_{jk}\}\sb{i,j,k\in\{1,2\}}$ would yield conditions for  
$\{\gt ^i_{jk}\}\sb{i,j,k\in\{1,2\}}$ similar to these in Proposition~\ref{3bi}.
Only this time $y$ would play the role of~$x$. Let us also observe that, if we
replaced $\s$ by $q\sp2\s$ in (\ref{comp}), then (\ref{comp}) would have no
solutions whatsoever as long as the negative powers of $x$ and $y$ are 
disallowed.
\hfill{$\Diamond$}\ere

\baselineskip15pt

\bre\label{azzar}\em
As in Section~2, equations (\ref{sol}) uniquely determine $\nabla\sp L$ from
$\nabla\sp R$, and vice-versa. Here, however, it happens despite the fact
that formula (\ref{azza}) is not fulfilled. (It follows from (\ref{za1}) that
$AZ(A\sp1)\subseteq xy A\sp1$.)
\hfill{$\Diamond$}\ere

\baselineskip16.5pt
Notice that, had we required the $\s$-compatibility condition to be satisfied on
 the whole bimodule $A\sp1$, then, in contrast with Proposition~\ref{3bi}, we would
obtain that the Christoffel symbols of a left connection have to fulfill the 
following equations:
\bea\label{3f}
&&F^1_{11} = x\left( y^3(-qf^1_{12} + f^1_{21})+yf^1_{11}-qy^2f^1_{22}\right)~, 
\nonumber \\
&& F^1_{12} = x^2(yf^1_{22} + y^2f^1_{12})~, \nonumber \\
&& F^1_{21} = x^2(q^2yf^1_{22} + y^2f^1_{21})~,  \nonumber \\
&& F^1_{22} = -q^2x^3f^1_{22}~, \nonumber \\  
&& F^2_{11} = f ^2_{11} +q y^4(f ^2_{22} -f ^1_{12} 
-\mbox{$\frac{3q}{1-q}$}f^1_{21})   
+ qy^2( f ^1_{11} - f ^2_{21} -q f ^2_{12}) ~, 
\nonumber \\ 
&& F ^2_{12} = x\left( q^2y^3(- f ^2_{22} +f ^1_{12}) + yf ^2_{12}  +qy^2f 
^1_{22}\right)  \nonumber \\ 
&& F^2_{21} =
x\left(qy^3( -f ^2_{22}  + q f ^1_{21})  +yf ^2_{21} +y^2f ^1_{22}\right)~, 
\nonumber \\ 
&& F^2_{22} = x^2(-q^2yf ^1_{22}+ y^2f ^2_{22})~, 
\eea
where  $f\sp i\sb{jk}\in Z(A),\, i,j,k\in\{1,2\}$.

The metric compatibility condition (\ref{g-comp-f}) can again be written in the form
of formula (\ref{compa}). This time, however, the centre of $A\sp1$ is 
non-trivial and it makes sense to ask what would happen if we imposed 
the metric compatibility condition only over $Z(A\sp1)$. It turns out that we have:

\bpr\label{exte}
Requiring that  the metric compatibility condition (\ref{g-comp-f}) be
satisfied only over $Z(A\sp1)$ is equivalent to demanding that it be fulfilled over
the whole $A\sp1$.
\epr 
{\it Proof.} Let $\mbox{\bf P}^{ij}$ denote the expression in the large parentheses
in (\ref{compa}). Furthermore, let us put 
$\mbox{\bf P}^{ij}=\th\sp k\mbox{\bf P}^{ij}\sb k$, $a_{1} = axy -bxy^3$, 
$a_2 = bx^2 y^2$, $\tilde{a}_{1}=a'xy -qb'xy^3$, $\tilde{a}_2=b'x^2y^2$, 
where $a, b, a', b'$ are arbitrary elements of $Z(A)$. 
Thanks to (\ref{za1}), substituting those terms
to (\ref{compa}) yields an equation that expresses the metric compatibility
condition over $Z(A\sp1)$: 
\bea
&&(axy -bxy^3)\th^k\mbox{\bf P}^{11}_k(a'xy -qb'xy^3) 
+ (axy -bxy^3)\th^k\mbox{\bf P}^{12}_kb'x^2y^2
\nonumber\\ &&
+ bx^2 y^2\th^k\mbox{\bf P}^{21}_k(a'xy -qb'xy^3)
+ bx^2 y^2\th^k\mbox{\bf P}^{21}_kb'x^2y^2=0\,
\label{zm}\eea
Commuting everything to the right of $\{\th\sp i\}\sb{i\in\{1,2\}}$ and taking
advantage of the fact that  $\{\th\sp i\}\sb{i\in\{1,2\}}$ is a basis, 
one can reduce (\ref{zm}) to:

{\samepage
\vspace*{-5mm}\baselineskip12pt\begin{eqnarray}
&& 
\Llp\Llp
xy\mbox{\bf P}^{11}_1 + (q-q^2)y^2\mbox{\bf P}^{11}_2
\Lrp a\nonumber \\ && +\Llp
-q^2x\mbox{\bf P}^{11}_1y^3 -(q^2-1)y\mbox{\bf P}^{11}_2y^3+x^2y^2
\mbox{\bf P}^{21}_1+(q^2-q)xy^3\mbox{\bf P}^{21}_2 
\Lrp b\Lrp 
(a'xy -qb'xy^3)\nonumber \\ &&  
+ \Llp\Llp
xy\mbox{\bf P}^{12}_1 + (q-q^2)y^2\mbox{\bf P}^{12}_2
\Lrp a\nonumber \\ && +\Llp
-q^2x\mbox{\bf P}^{12}_1y^3 -(q^2-1)y\mbox{\bf P}^{12}_2y^3+x^2y^2
\mbox{\bf P}^{22}_1 +(q^2-q)xy^3\mbox{\bf P}^{22}_2 
\Lrp b\Lrp b'x^2y^2 =0\nonumber\\
&& \ \nonumber \\
&&  \left(xy\mbox{\bf P}^{11}_2a-(qx
\mbox{\bf P}^{11}_2y^3 
-x^2y^2\mbox{\bf P}^{21}_2)b\right)(a'xy-qb'xy^3) \nonumber \\
&& + \left(xy\mbox{\bf P}^{12}_2a-(qx\mbox{\bf P}^{12}_2y^3 
-x^2y^2\mbox{\bf P}^{22}_2)b\right)b'x^2y^2=0
\end{eqnarray}\baselineskip16.5pt
Now, since $a,b,a',b'$ are arbitrary elements of $Z(A)$, the above two
equations boil down to:
}

\vspace*{-9mm}\baselineskip12pt\begin{eqnarray}
&&~\nonumber\\
&& xy\mbox{\bf P}^{11}_2xy=0 \nonumber \\
&& ~ \nonumber \\
&& xy\mbox{\bf P}^{11}_1xy +(q-q^2)y^2\mbox{\bf P}^{11}_2xy=0 \nonumber\\
&& ~ \nonumber \\
&&-qxy\mbox{\bf P}^{11}_2xy^3 +xy\mbox{\bf P}^{12}_2x^2y^2 =0\nonumber\\
&& ~ \nonumber \\
&& -qxy\mbox{\bf P}^{11}_1xy^3 +(1-q^2)y^2\mbox{\bf P}^{11}_2xy^3
+xy\mbox{\bf P}^{12}_1x^2y^2 +(q-q^2)y^2\mbox{\bf P}^{12}_2x^2y^2=0 \nonumber\\
&& ~ \nonumber \\
&& -qx\mbox{\bf P}^{11}_2xy^4 +x^2y^2\mbox{\bf P}^{21}_2xy=0 \nonumber \\
&& ~ \nonumber \\
&& -q^2x\mbox{\bf P}^{11}_1xy^4 -(q^2-1)y\mbox{\bf P}^{11}_2xy^4
+x^2y^2\mbox{\bf P}^{21}_1xy +(q^2-q)xy^3\mbox{\bf P}^{21}_2xy=0 \nonumber\\
&& ~ \nonumber \\
&&q^2x\mbox{\bf P}^{11}_2xy^4 -qx^2y^2\mbox{\bf P}^{21}_2xy^3
-qx\mbox{\bf P}^{12}_2x^2y^5 +x^2y^2\mbox{\bf P}^{22}_2x^2y^2=0 \nonumber\\
&& ~ \nonumber \\
&& x\mbox{\bf P}^{11}_1xy^6 +(1-q)y\mbox{\bf P}^{11}_2xy^6
-qx^2y^2\mbox{\bf P}^{21}_1xy^3 -(1-q^2)xy^3\mbox{\bf P}^{21}_2xy^3
-q^2x\mbox{\bf P}^{12}_1x^2y^5\nonumber\\
&& +(1-q^2)y\mbox{\bf P}^{12}_2x^2y^5
+x^2y^2\mbox{\bf P}^{22}_1x^2y^2 +(q^2-q)xy^3\mbox{\bf P}^{22}_2x^2y^2=0 
\eea\baselineskip16.5pt

\vspace*{-1mm}
Hence $\mbox{\bf P}^{ij}_k=0$ for any $i,j,k\in\{1,2\}$. Consequently,
$\mbox{\bf P}^{ij}=0$ for any $i,j\in\{1,2\}$, and (\ref{g-comp-f}) follows,
as claimed.
\epf

The same effect, though in a more trivial way, occurs in the setting of
Section~2.
\vspace*{-1mm}
\Section{Conclusions}\baselineskip16.5pt\parindent5mm
\vspace*{-1mm}

As we have demonstrated, restricting the domain of the $\s$-compatibility
condition $\nabla\sp L=\s\circ\nabla\sp R$ to $Z(A\sp1)$ yields
a theory of noncommutative linear connections
that coincides with the classical theory in the commutative case, does not
discriminate against the left or right structure of a bimodule 
and appears to be rich in the noncommutative set-up. 
It is easy to check, however, that this \mbox{$\sigma$-compatibility} equation
is not, in general, gauge covariant,
either when considering it over the whole bimodule~\cite{plane}, or when considering it only over the centre of
a bimodule. To provide a simple example, let us assume the setting of Section~4
and choose our gauge transformation to be
$\th\sp i\mapsto U\sp i\sb{\ j}\th\sp j,\, 
U=\mbox{\scriptsize $\pmatrix{1&x\cr 0&1}$}$.
Its action on $\nabla\sp L$ is given by the formula (cf.~p.547 and p.559 
in~\cite{conbook}) 
$N\mapsto dU.U\sp{-1}+UNU\sp{-1}$, where
$\nabla\sp L\th\sp i=N\sp i\sb{\ k}\ota\th\sp k,\, 
N=(N\sp i\sb{\ k})\in M\sb2(A\sp1)$.
It is clear that, although the pure gauge connections 
($\g\sp i\sb{jk}=0=\widetilde{\g}\sp i\sb{jk},\, i,j,k\in\{1,2\}$) satisfy
$\nabla\sp L=\s\circ\nabla\sp R$, the Christoffel symbol $\g\sp1\sb{12}$ of
the connection obtained by the action of $U$ on the left pure gauge connection
does not fulfill (\ref{3bif}). More precisely, 
$\g\sp1\sb{12}=1\in\hspace*{-3.75mm}\slash\hspace*{2mm} xA$. 
(This can be obtained from the equation
$(dU.U\sp{-1})\sp i\sb{\ k}\ota\th\sp k=\th\sp j\ota\th\sp k\g\sp i\sb{jk},
\, i\in\{1,2\}$.)

A different approach to the generalized permutation $\s$ has been suggested 
in~\cite{order}.  In the set-up of Section~5 in~\cite{order}, $\s$  is a
function of the connection. Consequently, to determine the space of connections
one needs to take into account all possible bimodule homomorphisms $\s$.
Allowing $\s$ to vary makes it possible
 to change it under the action of a gauge transformation.
If we take as a gauge transformation a bimodule automorphism $f:E\rightarrow E$
and define its action by the formula 
\be\label{bg}
\left(\nabla\sp L,\nabla\sp R,\s\right)\mapsto
\left((id\ota f\sp{-1})\!\circ\!\nabla\sp L\!\circ\! f,\;
(f\sp{-1}\ota id)\!\circ\!\nabla\sp R\!\circ\! f,\;
(id\ota f\sp{-1})\!\circ\!\s\!\circ\!(f\ota id)\right),
\ee
then the $\s$-compatibility condition $\nabla\sp L=\s\circ\nabla\sp R$ is gauge
covariant. Furthermore, since the centre of a bimodule is preserved by the
bimodule automorphisms, the $\s$-compatibility condition is also gauge
covariant when considered only over $Z(E)$. 
\footnote{We owe noticing this point to Michel Dubois-Violette.}
\ One should bear in mind, however, that, roughly speaking, the bimodule
automorphisms correspond to the `commutative sector' of the space of gauge
transformations.

A more radical point of view that might deserve a detailed investigation relies on
employing the metric compatibility condition (\ref{g-comp-f}) rather than the
equation $\nabla\sp L=\s\circ\nabla\sp R$  to relate $\nabla\sp L$ and
$\nabla\sp R$ uniquely and without (undesirable) restrictions on 
$(\nabla\sp L,\nabla\sp R)$. Clearly, for nondegenerate metrics, formulas 
(\ref{dg}) and (\ref{compa}) provide this kind of mutual dependance of $\nabla\sp L$ 
and $\nabla\sp R$. (In those cases, the nondegeneracy of a metric $g$ simply means
that $(g\sp{ij})$ is an invertible matrix.)

Finally, let us remark that it is plausible that in order to obtain a satisfactory
definition of a bimodule connection, one needs to use the language of 
quantum principal bundles, and, having understood and thoroughly worked out
the left-right relationship 
in this context (see Theorem~4.13 and Remark~4.14 in~\cite{mpps} and Appendix~B
in~\cite{bmq}), translate the solution(s) to bimodule terms.
\newpage
\appendix\parindent0pt\baselineskip16.125pt

\section{Appendix: Commutation Formulas}

Here we provide commutation formulas for the differential algebra $\Omega(A)$
that is defined with the help of (\ref{comm}).
Let $Q_{-1}=Q_{0}=0$ and, for $n>0$, let $Q_n=\sum_{k=1}^{n}q^{2(k-1)}$.
For any natural  $n\geq 0$, we have: 

\baselineskip12pt\vspace*{-4mm}\begin{eqnarray}
&& x^n\xi=q^{2n}\xi x^n,\;\;\; 
x^n\eta =q^n\eta x^n+(q^2-1)Q_n\xi x^{n-1}y,\;\;\; 
y^n\xi=q^{n}\xi y^n,\;\;\; 
y^n\eta =q^{2n}\eta y^n; ~~~~~~~~
\nonumber \\ 
&& ~ \nonumber \\
&& x^n\xi\ota\xi=q^{4n}\xi\ota\xi x^n ~, \nonumber \\
&& x^n\xi\ota\eta =q^{3n}\xi\ota\eta x^n+(q^2-1)q^{2n}Q_n\xi
\ota\xi x^{n-1}y~, \nonumber \\
&& x^n\eta \ota\xi=q^{3n}\eta \ota\xi x^n+(q^2-1)q^{2n-1}Q_n\xi
\ota\xi x^{n-1}y~, \nonumber \\
&& x^n\eta \ota\eta =q^{2n}\eta \ota\eta x^n+(q^2-1)q^{n}Q_n\eta 
\ota\xi x^{n-1}y  \nonumber \\
&& +(q^2-1)q^{n+1}Q_n\xi\ota\eta x^{n-1}y 
+q^2(q^2-1)^2 Q_nQ_{n-1}\xi \ota\xi x^{n-2}y^2~; \nonumber \\
&& ~ \nonumber \\
&& y^n\xi\ota\xi=q^{2n}\xi\ota\xi y^n~, ~~
y^n\xi\ota\eta =q^{3n}\xi\ota\eta y^n~, \nonumber \\ 
&& y^n\eta \ota\xi=q^{3n}\eta \ota\xi y^n~, 
~~ y^n\eta \ota\eta =q^{4n}\eta \ota\eta y^n~; \label{cf1}~~~~~~~~~~~~~~~~~~~~ 
\eea
\baselineskip16.125pt
One also has:
\bea
 && a\xi+b\eta:= 
 a_{pr}x^py^r\xi+b_{st} x^sy^t\eta  
\nonumber \\ 
&&=\xi\llp q^{2p + r} a_{pr}x^py^r
 + q^{2t}(q^2-1)Q_sb_{st}x^{s-1}y^{t+1}\lrp
+\eta \llp q^{2t+s}b_{st}x^sy^t\lrp  
\nonumber \\
&&= \xi\llp \llp q^{2p+r}a_{pr} +q^{2r-2}(q^2-1)Q_{p+1}
b_{p+1~r-1} \lrp x^py^r + q^{2p}a_{p0}x^p \lrp
+ \eta \llp q^{2t+s}b_{st}x^sy^t\lrp ~~~~~~~~~~~~~ 
\nonumber \\
&&=: \xi\tilde{a} + \eta \tilde{b}~, \label{cf2} 
\eea
where $a_{pr},b_{st}\in{\IC},\, p,r,s,t\in\{ 0,1,...\}$.

\section*{Acknowledgements} It is a pleasure to thank Michel Dubois-Violette
for helpful discussions.

\end{document}